\documentclass[a4paper,10pt,notitlepage,reqno]{article}
\usepackage[english]{babel}
\usepackage{amssymb,amsthm,bm} 
\usepackage{mathtools}
\usepackage{mathrsfs}
\usepackage{enumitem}
\usepackage{authblk}
\usepackage{cite}
\usepackage{xcolor}
\usepackage{tikz}
\usepackage{stackrel}
\usepackage{geometry}
\usepackage{hyperref}
\theoremstyle{plain} 
\newtheorem{theorem}{Theorem}[section]
\newtheorem*{theorem*}{Theorem}
\newtheorem{lemma}{Lemma}[section]
\newtheorem{proposition}{Proposition}[section]
\newtheorem{corollary}{Corollary}[section]
\newtheorem*{corollary*}{Corollary}

\theoremstyle{definition}
\newtheorem*{definition}{Definition}

\theoremstyle{remark}
\newtheorem{remark}{Remark}[section]

\newtheorem{example}{Example}
\numberwithin{equation}{section}

\usepackage{color}
\usepackage[normalem]{ulem}
\definecolor{DarkGreen}{rgb}{0,0.5,0.1} 

\newcommand\soutD{\bgroup\markoverwith
{\textcolor{DarkGreen}{\rule[.5ex]{2pt}{1pt}}}\ULon}

\newcommand{\Hm}[1]{\leavevmode{\marginpar{\tiny%
$\hbox to 0mm{\hspace*{-1.5mm}$\leftarrow$\hss}%
\vcenter{\vrule depth 0.1mm height 0.1mm width \the\marginparwidth}%
\hbox to
0mm{\hss$\rightarrow$\hspace*{-1.5mm}}$\\\relax\raggedright #1}}}

\definecolor{Darkgblue}{rgb}{0.3,0.3,0.5}

\newcommand{\la}{\lambda}

\newcommand{\ga}{\gamma}

\def\L{\mathcal{L}}

\newcommand{\erre}{\mathbb{R}}     

\newcommand{\rN}{\erre^{N}}

\newcommand{\D}{\rN \diff \{0\}}              

\newcommand{\CR}{C_c^{\infty}(\rN \diff \{0\})}   
\newcommand{\CD}{C_c^{\infty}(\D)}                

\newcommand{\pa}{\partial}

\newcommand{\diff}{\!\setminus\!}
\newcommand{\ov}{\overline}

\DeclarePairedDelimiter{\abs}{\lvert}{\rvert}

\DeclareMathOperator{\Div}{div}
\newcommand{\norma}[1]{\left\lVert#1\right\rVert}

\title{ A unified approach to \(L^p\) Hardy and Rellich-type inequalities in Euclidean and non-Euclidean settings}

\author[1]{Lorenzo D'Arca}

\affil[1]{Department of Mathematics “Guido Castelnuovo”, Sapienza University of Rome,\newline Piazzale Aldo Moro 5, Roma
00185, Italy; lorenzo.darca@uniroma1.it}

\date{}

\begin{document}

\maketitle

\begin{abstract}
We present a unified and concise method for establishing \(L^p\) Hardy and Rellich inequalities for a broad class of subelliptic operators of divergence type.  
The approach, based on a fundamental algebraic identity, provides explicit control on maximizing sequences and yields sharp constants in several significant cases.  
It applies beyond the Euclidean framework, covering the Heisenberg and Carnot group settings, and extends to a variety of subelliptic operators such as the Heisenberg–Greiner and Baouendi–Grushin operators.
\end{abstract}

\section{Introduction}
Hardy and Rellich inequalities are fundamental tools that naturally arise in a wide range of areas, including Analysis, Mathematical Physics, Spectral Theory, Geometry, and Quantum Mechanics. The sharp forms of these inequalities play a crucial role in establishing existence and nonexistence results for various partial differential equations. 
Due to their broad applicability, both inequalities have been extensively studied, giving rise to a vast body of literature (see, for instance, \cite{BE,OK,MMP,M,B,BT,CM} and the references therein). 

Hardy's inequality, introduced in \cite{H}, is one of the mathematical manifestations of the Uncertainty Principle. In its classical Euclidean form, it states that
\begin{equation}\label{Hardy Euc Lp Classica}
\left\lvert \frac{N - p}{p} \right\rvert^{p}
\int_{\mathbb{R}^N} \frac{|u|^p}{|x|^p} \, dx
\le
\int_{\mathbb{R}^N} |\nabla u|^p \, dx,
\quad \forall\, u \in C_c^\infty(\mathbb{R}^N \setminus \{0\}).
\end{equation}
A closely related version, involving only the radial component of the gradient, is given by
\begin{equation}\label{Hardy Radiale Euc Lp Classica}
\left\lvert\frac{N-p}{p} \right\rvert^p
\int_{\mathbb{R}^N}\frac{|u|^p}{|x|^p}\,dx
\leq
\int_{\mathbb{R}^N} \left\lvert \nabla u \cdot \frac{x}{|x|} \right\rvert^p \, dx.
\end{equation}
As \( |\nabla u \cdot x/|x|| \le |\nabla u| \), the inequality \eqref{Hardy Radiale Euc Lp Classica} is strictly stronger than \eqref{Hardy Euc Lp Classica}.
Nevertheless, both inequalities share the same optimal constant 
\( c_{N,p} = \left\lvert \frac{N - p}{p} \right\rvert^{p} \).
If \( N > p \), one has \( |x|^{-p} \in L^1_{\text{loc}}(\rN) \), 
and the validity of both \eqref{Hardy Euc Lp Classica} and \eqref{Hardy Radiale Euc Lp Classica} 
extends to all functions in the Sobolev space \(\mathcal{W}^{1,p}(\rN)\)\footnote{%
\(\mathcal{W}^{1,p}(\rN)\) denotes the Sobolev space of \(L^p(\rN)\) functions whose first weak derivatives also belong to \(L^p(\rN)\).}.
On the other hand, when \( N < p \), the inequalities hold for every \(u\) in the smaller domain
\(\mathcal{W}^{1,p}(\rN\diff\{0\})\coloneqq \ov{\CR}^{\norma{\cdot}_{\mathcal{W}^{1,p}(\rN)}}
\subsetneq \mathcal{W}^{1,p}(\rN)
\).
The dimension \(N=p\) is critical for the validity of the inequalities, as they cannot hold in this case with a non-zero constant on the left-hand side.
A systematic viewpoint on weighted Hardy inequalities was developed by Ghoussoub and Moradifam~\cite{GM0, GM1}, who introduced the notion of \emph{Bessel pairs} to characterize their validity.
Their approach offers a structural criterion, formulated in terms of an associated differential equation, that is closely related in spirit to the constructive viewpoint adopted in the present work.

A related higher-order inequality is Rellich's inequality, introduced in \cite{R}, which asserts that, for all 
\(u \in C_c^\infty(\rN \setminus \{0\})\),
\begin{equation}\label{Rellich Euc Lp Classica}
\left( \frac{N(p-1)(N-2p)}{p^2}\right)^p
\int_{\rN} \frac{\abs{u}^p}{\abs{x}^{2p}}\,dx
\leq
\int_{\rN} \abs{\Delta u}^p\,dx,
\qquad 2 \le p \le \frac{N}{2}.
\end{equation}
As in the Hardy case, the constant on the left-hand side of \eqref{Rellich Euc Lp Classica} is sharp.
When \(N > 2p\), \eqref{Rellich Euc Lp Classica} extends, by density, to all functions 
\(u \in W^{2,p}(\rN)\).
The dimension \(N = 2p\) plays the same critical role here as \(N = p\) does for Hardy's inequality.

\medskip
Given their central role in analysis and partial differential equations, considerable effort has been devoted to extending these inequalities to non-Euclidean settings. 
In the Heisenberg group \(\mathbb{H}^n\), the analogue of Hardy's inequality was established by Garofalo and Lanconelli in \cite{GL}, who proved that, for every 
\(u \in C_c^\infty(\mathbb{H}^n\setminus\{0\})\),
\begin{equation}
\label{Hardy Heisenberg Classica}
\frac{(Q-2)^2}{4}
\int_{\mathbb{H}^n} \frac{|u|^2}{\rho^2}\,|\nabla_{\mathbb{H}^n}\rho|^2 \, dzdt
\le 
\int_{\mathbb{H}^n} |\nabla_{\mathbb{H}^n} u|^2\, dzdt.
\end{equation}
The optimality of the constant was later proved by Goldstein and Zhang \cite{GZ}. 
Nonlinear extensions were obtained by Niu, Zhang and Wang \cite{NZW}, and further developments on groups of Heisenberg type are due to Danielli, Garofalo and Phuc \cite{DGP}. Improved Hardy inequalities in the presence of horizontal magnetic fields on \(\mathbb{H}^n\) 
have been recently derived in \cite{cassano2023horizontal}. 
Hardy-type inequalities on polarizable Carnot groups were established by D’Ambrosio \cite{DA0} and by Goldstein and Kombe \cite{GK}, while formulations on general Carnot groups can be found in the work of Jin and Shen \cite{JSW} and of Lian \cite{L}.

Higher-order analogues have also been studied in the framework of homogeneous groups. 
Rellich-type inequalities on stratified Lie groups were considered in \cite{CCR}, although the optimal constant was not determined. 
Sharp Rellich-type inequalities for the sub-Laplacian in Carnot groups were subsequently obtained by Kombe \cite{K} and by Lian \cite{L}, with further results in \cite{DA1,JH,RS,NLN,N}.
Relevant applications of Hardy and Rellich inequalities to linear and nonlinear equations on Carnot groups appear in several works, including \cite{GL}, \cite{DA1}, \cite{AA}, \cite{GK}, \cite{GK0}, \cite{GZ}, \cite{JH}, \cite{JZ}, \cite{K} and related references.

\medskip
The proofs developed in this work rely on the non-negativity of a suitable \(L^2\) norm, a technique that lies at the core of many well-known arguments for Hardy and Rellich inequalities in the quadratic setting. 
In particular, as explained in \cite{C}, several classical inequalities follow from the elementary identity  
\[
0 \le \|f-g\|_{L^2}^2 
   = \|f\|_{L^2}^2 + \|g\|_{L^2}^2 - 2(f,g)_{L^2}.
\]

To treat the general \(L^p\) case, we employ a nonlinear extension of this argument. 
More precisely, there exists a positive weight \(w(p,f,g)\geq 0\) such that  
\[
0 \le \|w(p,f,g)(f-g)\|_{L^2}^2
   = \|f\|_{L^p}^p + (p-1)\|g\|_{L^p}^p 
     - p\bigl(|g|^{p-2}g,f\bigr)_{L^2}.
\]
The explicit expression and basic properties of \(w(p,f,g)\) will be recalled in Section~\ref{Sezione disuguaglianze Lp}.  
In particular, by definition one has \(w \ge 0\) and \(w = 0\) if and only if \(f = g = 0\), so that equality in the above relations occurs precisely when \(f=g\).

Identities of this type have appeared in related contexts, for instance in the works of Ioku, Ishiwata and Ozawa~\cite{IIO} in the Euclidean setting, as well as in the work of Ruzhansky and Suragan~\cite{RS} for homogeneous groups.

The approach developed here allows us to obtain explicit formulas for the maximizing functions and to prove the sharpness of the associated constants in several significant cases.
These results provide a unified framework that both encompasses and refines several known Hardy and Rellich inequalities, and establish general sufficient conditions for the validity and sharpness of weighted inequalities associated with subelliptic operators that do not necessarily arise from a homogeneous group structure.

\medskip
Before stating our main results, we outline the analytical framework in which they are developed.
We work in \(\mathbb{R}^N\) and consider \(h \le N\) vector fields of the form
\[
X_i(x) = \sum_{j=1}^{N} \sigma_{i,j}(x)\,\frac{\partial}{\partial x_j}, 
\qquad i = 1, \ldots, h,
\]
where each coefficient \(\sigma_{i,j}\) is continuous on \(\mathbb{R}^N\) and its partial derivatives 
\(\frac{\partial}{\partial x_j}\sigma_{i,j}\) are continuous as well.  
Let
\[
\sigma(x) = \big(\sigma_{i,j}(x)\big)_{i=1,\ldots,h}^{j=1,\ldots,N}
\]
denote the associated \(h\times N\) matrix.  
We define the first-order differential operator
\begin{equation}\label{def:nablaL}
\nabla_{\L} \coloneqq (X_1,\ldots,X_h) = \sigma(x)\nabla,
\end{equation}
and denote by \(X_i^*\) the formal adjoint of \(X_i\) with respect to the Lebesgue measure
\[
X_i^* = - \sum_{j=1}^{N} \frac{\partial}{\partial x_j}\big(\sigma_{i,j}(x)\,\cdot\,\big),
\qquad i = 1,\ldots,h.
\]

The symmetric matrix
\[
A(x) \coloneqq \sigma(x)^{T}\sigma(x)
\]
naturally arises in this setting, and allows us to introduce the second-order operator
\begin{equation}\label{def:L}
\L \coloneqq \Div(A(x)\nabla)
= -\sum_{i=1}^{h} X_i^*X_i
= -\,\nabla_{\L}^*\cdot\nabla_{\L}.
\end{equation}
For \(p>1\), we also consider the nonlinear operator
\begin{equation}\label{def:Lp}
\L_p(u) \coloneqq \Div\big(|\nabla_{\L}u|^{p-2} A(x)\nabla u\big)
= -\,\nabla_{\L}^* \cdot \big(|\nabla_{\L}u|^{p-2}\nabla_{\L}u\big).
\end{equation}

\begin{example}[Euclidean Laplacian]
If \(X_i = \frac{\partial}{\partial x_i}\) for \(i=1,\ldots,N\), 
then \(\nabla_{\L}\) coincides with the standard gradient \(\nabla\), 
\(A(x)\) is the identity matrix on \(\mathbb{R}^N\), 
and \(\L\) reduces to the classical Laplacian \(\Delta\). 
In this setting, \(\L_p\) corresponds to the standard \(p\)-Laplacian.
\end{example}

Now we are in position to state our main results. 
In the following, we assume that \(d\colon \D\to (0,+\infty)\) is a smooth, positive, non-constant function.

\begin{theorem}[Weighted Hardy inequality]
\label{INTRO: SV: TH Hardy Pesata}
Let \(p \geq 2\) and \(\beta \in \mathbb{R}\) be such that 
\[
\L_p d = \frac{(1 - \beta)(p - 1)\,|\nabla_{\L} d|^p}{d}
\quad \text{in } \D.
\]
Then, for every \(\theta \in \mathbb{R}\) and \(u \in \CD\),
the following inequality holds:
\begin{equation}
\label{INTRO: SV: Hardy Pesata}
\left| \frac{p(\theta - 1) + \beta(p - 1)}{p} \right|^p
\int_{\rN} \frac{|u|^p}{d^{p\theta}}\,|\nabla_{\L} d|^p\,dx
\leq
\int_{\rN} 
\left|\nabla_{\L} u \cdot 
\frac{\nabla_{\L} d}{|\nabla_{\L} d|}\right|^p
\frac{1}{d^{p(\theta - 1)}}\,dx
\leq
\int_{\rN} 
\frac{|\nabla_{\L} u|^p}{d^{p(\theta - 1)}}\,dx.
\end{equation}
Moreover, if there exists a positive constant \(\gamma_p\) such that
\[
\int_{\{d = r\}} 
\frac{|\nabla_{\L} d|^p}{|\nabla d|}\,d\mathcal{H}_{N-1}
= \gamma_p\,r^{(p - 1)(1 - \beta)}
\qquad \text{for a.e. } r > 0,
\]
then the constant on the left-hand side is sharp, and the (non-admissible) extremal function is, up to a multiplicative constant,
\[
u(x) = d(x)^{\frac{p(\theta - 1) + \beta(p - 1)}{p}}.
\]
\end{theorem}

\begin{theorem}[Weighted Rellich inequality]
\label{INTRO: SV: TH Rellich Pesata}
Let \(p \ge 2\), \(\beta \in \mathbb{R}\), and \(\theta \in \mathbb{R}\) be such that
\[
\L d = \frac{(1 - \beta)\,|\nabla_{\L} d|^2}{d}
\quad \text{in } \D,
\qquad
d^{-p\theta}\,|\nabla_{\L} d|^{-2(p-1)} \in L^1_{\mathrm{loc}}(\D),
\]
and
\[
(2 - \beta - p\theta - 2p)\,(p\theta + 2p - 2 - \beta(p - 1)) \ge 0.
\]
Then, for every \(u \in \CD\), the following inequality holds:
\begin{equation}
\label{INTRO: SV: Rellich Pesata}
\left(
\frac{(2 - \beta - p\theta - 2p)\,(p\theta + 2p - 2 - \beta(p - 1))}{p^2}
\right)^p
\int_{\rN} \frac{|u|^p}{d^{p(\theta + 2)}}\,|\nabla_{\L} d|^2\,dx
\le
\int_{\rN} \frac{|\L u|^p}{d^{p\theta}}\,\frac{1}{|\nabla_{\L} d|^{2(p-1)}}\,dx.
\end{equation}
Moreover, if there exists a positive constant \(\lambda\) such that
\[
\int_{\{d = r\}} 
\frac{|\nabla_{\L} d|^2}{|\nabla d|}\,d\mathcal{H}_{N-1}
= \lambda\,r^{1 - \beta}
\qquad \text{for a.e. } r > 0,
\]
then the constant on the left-hand side of \eqref{INTRO: SV: Rellich Pesata} is sharp, and the (non-admissible) extremal function is, up to a multiplicative constant,
\[
u(x) = d(x)^{\frac{\beta + p\theta + 2p - 2}{p}}.
\]
\end{theorem}

\begin{remark}
The assumptions on \(\L_p d\) and \(\L_2 d\) in 
Theorems~\ref{INTRO: SV: TH Hardy Pesata} and~\ref{INTRO: SV: TH Rellich Pesata} 
can be conveniently verified starting from conditions of the form
\[
\begin{cases}
\L_p(d^{\beta}) = 0 & \text{in } \D, \quad \text{if } \beta \neq 0,\\[4pt]
\L_p(-\ln d) = 0 & \text{in } \D, \quad \text{if } \beta = 0,
\end{cases}
\]
see Lemmas~\ref{SV: p-Laplaciano di d^beta} and~\ref{SV: p-Laplaciano di ln}. 
Furthermore, the integral condition
\[
\int_{\{d = r\}} \frac{|\nabla_{\L} d|^p}{|\nabla d|}\,d\mathcal{H}_{N-1}
= \gamma_p\, r^{(p-1)(1-\beta)}, \qquad \text{for a.e. } r>0,
\]
which appears in the sharpness statements of both 
Theorems~\ref{INTRO: SV: TH Hardy Pesata} and~\ref{INTRO: SV: TH Rellich Pesata},
is fulfilled, for instance, when \(d\) and \(\nabla_{\L}\) satisfy suitable homogeneity relations 
with respect to a family of dilations, as discussed in the introductory part of 
Section~\ref{sezione: applicazioni}.
\end{remark}

In contrast to the works of Ruzhansky and Suragan~\cite{RS} and Nguyen et al.~\cite{NLN}, which are formulated in terms of the Euler operator associated with a family of dilations on homogeneous groups, 
the inequalities established here are expressed purely in terms of horizontal derivatives.
This allows us to treat subelliptic structures where no natural Euler-type operator is available, or where the underlying space does not possess a homogeneous group structure.

\subsection{Organization of the Paper}
The paper is organized as follows.
In Section~\ref{sec:tools} we collect the algebraic identities and analytic inequalities that constitute the technical core of the work.
Section~\ref{SEZIONE: GENERALE} contains the proof of the main results, Theorems~\ref{INTRO: SV: TH Hardy Pesata} and~\ref{INTRO: SV: TH Rellich Pesata}.
For clarity, the argument is divided into two parts: the validity of the inequalities is established in Section~\ref{sezione: validita}, while Section~\ref{section: sharpness constant} is devoted to the analysis of the sharpness of the corresponding constants.
Finally, Section~\ref{sezione: applicazioni} illustrates several model cases in which the general theorems apply, including the Euclidean setting, the Heisenberg–Greiner operator, the Baouendi–Grushin operator, and the sub-Laplacian on Carnot groups.

\section{Auxiliary tools}\label{sec:tools}
\subsection{Algebraic Identities}\label{Sezione disuguaglianze Lp}
The core of our approach relies on a purely algebraic identity. 
Although similar formulas have appeared in related contexts 
(see, e.g., \cite{IIO, RS}), 
we recall here a self-contained version adapted to our notation, 
as it constitutes the starting point for the subsequent functional inequalities.

\begin{proposition}[Fundamental \(L^p\) identity]\label{Proposizione identità Lp}
Let \(D \subseteq \mathbb{R}^N\) be an open set, and let \(p \ge 2\).
For all \(f,g \in L^p(D)\), the following identity holds:
\begin{equation}\label{Idenità L^p}
\|w(p,f,g)(f-g)\|_{L^2(D)}^2
= \|f\|_{L^p(D)}^p + (p-1)\|g\|_{L^p(D)}^p
- p\,(\,|g|^{p-2}g,\,f\,)_{L^2(D)},
\end{equation}
where the weight \(w(p,f,g)\) is defined by
\begin{equation}\label{Definizione peso w}
w(p,f,g)^2
\coloneqq p(p-1)\int_0^1 s\,|\,s g + (1-s)f\,|^{p-2}\,ds.
\end{equation}
\end{proposition}

\begin{proof}
Starting from the term 
\[
|f|^p + (p-1)|g|^p - p|g|^{p-2} g f,
\]
one may add and subtract \(p|f|^p\) to obtain
\[
(p-1)(|g|^p - |f|^p) - p f \big(|g|^{p-2}g - |f|^{p-2}f\big).
\]
Noting that
\[
|g|^p - |f|^p
= \int_0^1 \frac{\partial}{\partial s}|s g+(1-s)f|^p\,ds,
\quad
|g|^{p-2}g - |f|^{p-2}f
= \int_0^1 \frac{\partial}{\partial s}\!\big(|s g+(1-s)f|^{p-2}(s g+(1-s)f)\big)\,ds,
\]
we arrive at
\[
|f|^p + (p-1)|g|^p - p|g|^{p-2} g f
= p(p-1)(g-f)^2 \int_0^1 s\,|s g+(1-s)f|^{p-2}\,ds.
\]
Integration over \(D\) yields the asserted identity.
\end{proof}

\begin{remark}
By construction \(w(p,f,g) \ge 0\), and \(w(p,f,g) = 0\) if and only if \(f=g=0\).
In particular, equality in Proposition~\ref{Proposizione identità Lp} occurs if and only if \(f=g\).
\end{remark}

Another elementary but useful estimate is contained in the following lemma.

\begin{lemma}\label{Stima di |a+b|^p}
For every \(p \ge 2\) there exists a constant \(c_p > 0\) such that, for all \(a,b \in \mathbb{R}\),
\[
|a+b|^p \le |a|^p + c_p \bigl(|a|^{p-1}|b| + |b|^p\bigr).
\]
\end{lemma}

\begin{proof}
The claim follows by setting \(x = a/b\) (for \(b \ne 0\)) and observing that
\[
f(x) = \frac{|x+1|^p - |x|^p}{|x|^{p-1}+1}
= p \int_0^1 \frac{|s+x|^{p-2}(s+x)}{|x|^{p-1}+1}\,ds
\]
is continuous and bounded on \(\mathbb{R}\).
\end{proof}

A double application of Lemma~\ref{Stima di |a+b|^p} yields the following corollary.

\begin{corollary}\label{Stima di |a+b+c|^p}
For every \(p \ge 2\) there exists a constant \(c_p > 0\) such that, for all \(a,b,c \in \mathbb{R}\),
\[
|a+b+c|^p
\le |a|^p
+ c_p |a|^{p-1}(|b|+|c|)
+ c_p |b|^p
+ c_p^2 |b|^{p-1}|c|
+ c_p^2 |c|^p.
\]
\end{corollary}

\subsection{Smooth cut-off functions}

We shall frequently employ a family of smooth cut-off functions that will be used in the construction of maximizing sequences and in the analysis of the sharpness of our inequalities. 
Given \(\varepsilon > 0\), let \(g_\varepsilon \in C_c^\infty(\mathbb{R}^+)\) satisfy
\begin{equation}\label{Funzioni Cut-Off}
g_\varepsilon(r)=
\begin{cases}
0, & 0\le r\le \varepsilon \text{ or } r\ge \dfrac{1}{\varepsilon},\\[6pt]
1, & 2\varepsilon \le r \le \dfrac{1}{2\varepsilon},
\end{cases}
\qquad 0\le g_\varepsilon(r)\le 1 \text{ for all } r\ge0,
\end{equation}
and
\[
\begin{cases}
|g'_\varepsilon(r)|\le \dfrac{c}{\varepsilon},\quad |g''_\varepsilon(r)|\le \dfrac{c}{\varepsilon^2}, & \varepsilon \le r \le 2\varepsilon,\\[6pt]
|g'_\varepsilon(r)|\le c\varepsilon,\quad |g''_\varepsilon(r)|\le c\varepsilon^2, & \dfrac{1}{2\varepsilon}\le r \le \dfrac{1}{\varepsilon},
\end{cases}
\]
for some constant \(c>0\). 
A straightforward computation shows that
\begin{equation}\label{Integrali g_epsilon}
\begin{aligned}
&\int_{0}^{+\infty}r^{-1}|g_\varepsilon(r)|^p\,dr = -\ln(4\varepsilon^2)+\mathcal{O}(1),\\
&\int_{0}^{+\infty}|g_\varepsilon(r)|^{p-1}|g'_\varepsilon(r)|\,dr= \mathcal{O}(1),\\
&\int_{0}^{+\infty}r^{p-1}|g'_\varepsilon(r)|^p\,dr = \mathcal{O}(1),\\
&\int_{0}^{+\infty}r|g_\varepsilon(r)|^{p-1}|g''_\varepsilon(r)|\,dr = \mathcal{O}(1),\\
&\int_{0}^{+\infty}r^p|g'_\varepsilon(r)|^{p-1}|g''_\varepsilon(r)|\,dr = \mathcal{O}(1),\\
&\int_{0}^{+\infty}r^{2p-1}|g''_\varepsilon(r)|^p\,dr = \mathcal{O}(1).
\end{aligned}
\end{equation}

\section{Proof of the main results}\label{SEZIONE: GENERALE}

We divide the proofs of Theorems~\ref{INTRO: SV: TH Hardy Pesata} and~\ref{INTRO: SV: TH Rellich Pesata} into two parts.
The first one establishes the validity of the inequalities, while the second one concerns the sharpness of the corresponding constants.

\subsection{Validity of the weighted Hardy and Rellich inequalities}\label{sezione: validita}
In this subsection we establish the validity of the inequalities stated in 
Theorems~\ref{INTRO: SV: TH Hardy Pesata} and~\ref{INTRO: SV: TH Rellich Pesata}.
We first deal with the Hardy inequality, and then with its higher–order counterpart.

\begin{proposition}[Validity of the weighted Hardy inequality]
\label{Proof: Hardy Validity}
Let \(p \ge 2\), \(\beta, \theta \in \mathbb{R}\), and assume that 
\[
\L_p d = \frac{(1 - \beta)(p - 1)\,|\nabla_{\L} d|^p}{d}
\quad \text{in } \D.
\]
Then the inequality stated in 
Theorem~\ref{INTRO: SV: TH Hardy Pesata} holds for every 
\(u \in C_c^\infty(\D)\).
\end{proposition}

\begin{proof}
We apply identity~\eqref{Idenità L^p} to
\[
f = \Bigl(\nabla_{\L} u \cdot \frac{\nabla_{\L} d}{|\nabla_{\L} d|}\Bigr)\, d^{1-\theta},
\qquad
g = \alpha\, u\,|\nabla_{\L} d|\, d^{-\theta},
\quad \alpha\in\mathbb{R}.
\]
We obtain
\[
\begin{aligned}
0 \le {} 
\| w(p,f,g)\,(f-g) \|_{L^2}^2
= {} &
\int_{\mathbb{R}^N}
\Bigl|\nabla_{\L} u \cdot \frac{\nabla_{\L} d}{|\nabla_{\L} d|}\Bigr|^p \frac{1}{d^{p(\theta-1)}}\,dx
+
(p-1)|\alpha|^p
\int_{\mathbb{R}^N}
\frac{|u|^p}{d^{p\theta}}\,
|\nabla_{\L} d|^p\,dx
\\[4pt]
& \quad
-
p|\alpha|^{p-2}\alpha
\int_{\mathbb{R}^N}
\frac{|u|^{p-2}u}{d^{p\theta-1}}\,
\bigl(\nabla_{\L} u \cdot \nabla_{\L} d\bigr)\,
|\nabla_{\L} d|^{p-2}\,dx.
\end{aligned}
\]
Using \(\nabla_{\L}|u|^p = p|u|^{p-2}u\,\nabla_{\L}u\) and integrating by parts gives
\[
\begin{aligned}
\int_{\mathbb{R}^N}
\frac{|u|^{p-2}u}{d^{p\theta-1}}\,
\bigl(\nabla_{\L} u \cdot \nabla_{\L} d\bigr)\,
|\nabla_{\L} d|^{p-2}\,dx
&=
-\frac{1}{p}
\int_{\mathbb{R}^N}
\frac{|u|^p}{d^{p\theta-1}}\,
\L_p d\,dx
\\[4pt]
&\quad
+\frac{p\theta-1}{p}
\int_{\mathbb{R}^N}
\frac{|u|^p}{d^{p\theta}}\,
|\nabla_{\L} d|^p\,dx.
\end{aligned}
\]
Since
\(
\L_p d = \frac{(1-\beta)(p-1)\,|\nabla_{\L} d|^p}{d}
\)
in \(\D\),
we deduce
\[
\int_{\mathbb{R}^N}
\frac{|u|^{p-2}u}{d^{p\theta-1}}\,
\bigl(\nabla_{\L} u \cdot \nabla_{\L} d\bigr)\,
|\nabla_{\L} d|^{p-2}\,dx
=
\frac{p(\theta-1)+\beta(p-1)}{p}
\int_{\mathbb{R}^N}
\frac{|u|^p}{d^{p\theta}}\,
|\nabla_{\L} d|^p\,dx.
\]
Substituting yields
\[
\int_{\mathbb{R}^N}
\Bigl|\nabla_{\L} u \cdot \frac{\nabla_{\L} d}{|\nabla_{\L} d|}\Bigr|^p\,\frac{1}{d^{p(\theta-1)}}\,dx
\ge
-\Bigl[(p-1)|\alpha|^p
-
\bigl(p(\theta-1)+\beta(p-1)\bigr)
|\alpha|^{p-2}\alpha\Bigr]
\int_{\mathbb{R}^N}
\frac{|u|^p}{d^{p\theta}}\,
|\nabla_{\L} d|^p\,dx.
\]
Optimizing in \(\alpha\) gives
\(\displaystyle \alpha = \frac{p(\theta-1)+\beta(p-1)}{p}\),
whence
\[
\Bigl|\frac{p(\theta-1)+\beta(p-1)}{p}\Bigr|^p
\int_{\mathbb{R}^N}
\frac{|u|^p}{d^{p\theta}}\,
|\nabla_{\L} d|^p\,dx
\le
\int_{\mathbb{R}^N}
\Bigl|\nabla_{\L} u \cdot \frac{\nabla_{\L} d}{|\nabla_{\L} d|}\Bigr|^p \frac{1}{d^{p(\theta-1)}}\,dx.
\]
Finally, since
\(\bigl|\nabla_{\L} u \cdot \tfrac{\nabla_{\L} d}{|\nabla_{\L} d|}\bigr|\le |\nabla_{\L} u|\),
the second inequality in
\eqref{INTRO: SV: Hardy Pesata}
follows immediately.
\end{proof}

\begin{remark}
\label{Final Remark Hardy}
Equality in the first inequality of
\eqref{INTRO: SV: Hardy Pesata}
holds if and only if
\[
\Bigl(\nabla_{\L} u \cdot \frac{\nabla_{\L} d}{|\nabla_{\L} d|}\Bigr)\, d^{1-\theta}
=
\frac{p(\theta-1)+\beta(p-1)}{p}\,
u\,|\nabla_{\L} d|\, d^{-\theta}
\quad \text{in } \D.
\]
In particular, all \(d\)-radial solutions of this equation are of the form
\[
u(x) = \lambda\, d(x)^{\frac{p(\theta-1)+\beta(p-1)}{p}},
\qquad \lambda\in\mathbb{R}\setminus\{0\},
\]
and represent, up to scaling, the non-admissible functions
that make both equalities hold.
\end{remark}

Before turning to the Rellich inequality, we record a useful consequence of Proposition \ref{Proof: Hardy Validity}, which will play a central role in its proof.

\begin{lemma}
\label{SV: Corollario HardyXRellich}
Let \(p \ge 2\) and \(\theta, \beta \in \mathbb{R}\).
Assume that
\[
\L d \;=\; \frac{(1 - \beta)\,|\nabla_{\L} d|^2}{d}
\quad \text{in } \D.
\]
Then, for every \(u \in \mathcal{C}_c^\infty(\D)\),
\begin{equation}
\label{SV: EQ: In Corollario HardyXRellich}
\int_{\mathbb{R}^N}
\frac{|u|^{p-2}\,|\nabla_{\L} u|^2}{d^{p\theta + 2p - 2}}\,dx
\;\ge\;
\frac{(p\theta + 2p - 2 + \beta)^2}{p^2}
\int_{\mathbb{R}^N}
\frac{|u|^p}{d^{p(\theta + 2)}}\,|\nabla_{\L} d|^2\,dx.
\end{equation}
\end{lemma}

\begin{proof}
From Proposition~\ref{Proof: Hardy Validity}, applied for \(p=2\) and with the parameter
\(\theta' = \tfrac{p}{2}(\theta+2)\), we obtain
\[
\int_{\mathbb{R}^N}
\frac{|\nabla_{\L} u|^2}{d^{p\theta + 2p - 2}}\,dx
\;\ge\;
\frac{(p\theta + 2p - 2 + \beta)^2}{4}
\int_{\mathbb{R}^N}
\frac{|u|^2}{d^{p(\theta + 2)}}\,|\nabla_{\L} d|^2\,dx.
\]
Finally, let \(\varepsilon>0\) and set
\[
u_\varepsilon \coloneqq \bigl(u^2+\varepsilon^2\bigr)^{\frac{p}{4}} - \varepsilon^{\frac{p}{2}}
\;\in\; \mathcal{C}_c^\infty(\D).
\]
Applying the preceding inequality to \(u_\varepsilon\) and letting \(\varepsilon \to 0^+\),
we obtain \eqref{SV: EQ: In Corollario HardyXRellich}.
\end{proof}

\begin{proposition}[Validity of the weighted Rellich inequality]
\label{Proof: Rellich Validity}
Assume the hypotheses of 
Theorem~\ref{INTRO: SV: TH Rellich Pesata}.
Then the inequality stated therein holds for every 
\(u \in C_c^\infty(\D)\).
\end{proposition}

\begin{proof}
We invoke identity~\eqref{Idenità L^p} with the choice
\[
f = \frac{\L u}{d^{\theta}}\,|\nabla_{\L} d|^{-2\frac{p-1}{p}},
\qquad
g = \alpha\,\frac{u}{d^{\theta + 2}}\,|\nabla_{\L} d|^{\frac{2}{p}},
\]
where \(\alpha \in \mathbb{R}\) is a free parameter to be optimized later.  
Integrating by parts, we find
\[
\int_{\mathbb{R}^N}
\frac{|\L u|^p}{d^{p\theta}\,|\nabla_{\L} d|^{2(p-1)}}\,dx
\ge
-(p-1)|\alpha|^p
\int_{\mathbb{R}^N}
\frac{|u|^p}{d^{p(\theta + 2)}}\,|\nabla_{\L} d|^2\,dx
+
p|\alpha|^{p-2}\alpha\,\mathcal{I},
\]
where
\[
\mathcal{I}
= \int_{\mathbb{R}^N}
\frac{|u|^{p-2}u\,\L u}{d^{p\theta + 2p - 2}}\,dx.
\]

A direct computation yields
\[
\mathcal{I}
= - (p-1)\!\int_{\mathbb{R}^N}
\frac{|u|^{p-2}|\nabla_{\L} u|^2}{d^{p\theta + 2p - 2}}\,dx
-
\frac{(p\theta + 2p - 2)(2 - \beta - p\theta - 2p)}{p}
\!\int_{\mathbb{R}^N}
\frac{|u|^p}{d^{p(\theta + 2)}}\,|\nabla_{\L} d|^2\,dx.
\]

Substituting this identity, we obtain
\[
\begin{aligned}
\int_{\mathbb{R}^N}
\frac{|\L u|^p}{d^{p\theta}\,|\nabla_{\L} d|^{2(p-1)}}\,dx
\ge&
-\Bigl\{
(p-1)|\alpha|^p
+\!|\alpha|^{p-2}\alpha(p\theta + 2p - 2)(2 - \beta - p\theta - 2p)
\Bigr\}
\!\int_{\mathbb{R}^N}\!\frac{|u|^p}{d^{p(\theta + 2)}}|\nabla_{\L} d|^2\,dx
\\[-2pt]
&-\,|\alpha|^{p-2}\alpha\,p(p-1)
\!\int_{\mathbb{R}^N}\!\frac{|u|^{p-2}|\nabla_{\L} u|^2}{d^{p\theta + 2p - 2}}\,dx.
\end{aligned}
\]

We restrict to \(\alpha \le 0\) so that 
Lemma~\ref{SV: Corollario HardyXRellich} applies to the last term 
with the correct sign.
The parameter condition 
\((2 - \beta - p\theta - 2p)\,(p\theta + 2p - 2 - \beta(p - 1)) \ge 0\)
guarantees that the optimizer obtained below indeed satisfies \(\alpha \le 0\).
Using Lemma~\ref{SV: Corollario HardyXRellich}, we obtain
\[
\begin{aligned}
\int_{\mathbb{R}^N}
\frac{|\L u|^p}{d^{p\theta}\,|\nabla_{\L} d|^{2(p-1)}}\,dx
\ge&
-\Biggl\{
(p-1)|\alpha|^p
+\!|\alpha|^{p-2}\alpha
\biggl[
(p\theta + 2p - 2)(2 - \beta - p\theta - 2p)
\\[-2pt]
&\qquad\qquad\qquad
+\frac{(p-1)(2 - \beta - p(\theta + 2))^2}{p}
\biggr]
\!\Biggr\}
\int_{\mathbb{R}^N}
\frac{|u|^p}{d^{p(\theta + 2)}}\,|\nabla_{\L} d|^2\,dx.
\end{aligned}
\]

Maximizing the term inside the braces with respect to \(\alpha\), we find that the optimal choice is
\[
\alpha_* = 
\frac{(p\theta + 2p - 2 + \beta)\,(p\theta + 2p - 2 - \beta(p - 1))}{p^2}.
\]
Substituting this value into the previous inequality yields 
\eqref{INTRO: SV: Rellich Pesata}.
\end{proof}

\begin{remark}
\label{Final remark Rellich}
Equality in~\eqref{INTRO: SV: Rellich Pesata} requires that \(f=g\) in~\eqref{Idenità L^p}, 
which yields the differential equation
\[
\frac{\L u}{d^{\theta}}\,|\nabla_{\L} d|^{-2\frac{p-1}{p}}
=
\frac{(p\theta + 2p - 2 + \beta)\,(p\theta + 2p - 2 - \beta(p - 1))}{p^2}
\,\frac{u}{d^{\theta + 2}}\,|\nabla_{\L} d|^{\frac{2}{p}}
\quad \text{in } \D.
\]
Its \(d\)-radial solutions are 
\[
u(x)=\lambda\,d(x)^{\frac{\beta+p\theta+2p-2}{p}},
\qquad \lambda\neq0,
\]
and these profiles also satisfy the equality condition in~\eqref{SV: EQ: In Corollario HardyXRellich}.
Hence they represent, up to scaling, the non-admissible functions 
for which equality holds in~\eqref{INTRO: SV: Rellich Pesata}.
\end{remark}

\subsection{Sharpness of the weighted Hardy and Rellich inequalities}\label{section: sharpness constant}
In this subsection we establish the sharpness of the constants in the weighted Hardy and Rellich inequalities stated in Theorems~\ref{INTRO: SV: TH Hardy Pesata} and~\ref{INTRO: SV: TH Rellich Pesata}. The argument is constructive: it relies on the smooth cut-off families introduced in Section~\ref{sec:tools} and on the explicit extremal profiles identified in Section~\ref{sezione: validita}.

\begin{proposition}[Sharpness of the weighted Hardy constant]\label{PROP: Hardy Sharpness}
Under the assumptions of Theorem~\ref{INTRO: SV: TH Hardy Pesata}, 
the constant in~\eqref{INTRO: SV: Hardy Pesata} is sharp.
\end{proposition}

\begin{proof}
Let 
\[
u_\varepsilon(x)
= d(x)^{\frac{p(\theta-1)+\beta(p-1)}{p}}\,g_\varepsilon(d(x)),
\qquad \varepsilon>0,
\]
where \(g_\varepsilon\) is the smooth cut-off function introduced in~\eqref{Funzioni Cut-Off}.  
Since \(u_\varepsilon\) depends only on \(d(x)\), we have
\[
\nabla_{\L}u_\varepsilon
= d^{\alpha_*-1}\bigl(\alpha_* g_\varepsilon + d g_\varepsilon'\bigr)\nabla_{\L}d,
\qquad 
\alpha_*=\frac{p(\theta-1)+\beta(p-1)}{p}.
\]
Applying Lemma~\ref{Stima di |a+b|^p} and invoking the coarea formula 
(see \cite[Theorem~3.2.12, p.~249]{Fe}) together with
\[
\int_{\{d=r\}}\frac{|\nabla_{\L} d|^p}{|\nabla d|}\,d\mathcal{H}_{N-1}
=\gamma_p\,r^{(p-1)(1-\beta)} \qquad \text{for a.e. } r>0,
\]
we obtain, for suitable constants \(c_1,c_2>0\),
\[
\begin{split}
\int_{\mathbb{R}^N}
\Bigl|\nabla_{\L} u_\varepsilon \cdot 
\frac{\nabla_{\L} d}{|\nabla_{\L} d|}\Bigr|^{p}
\frac{1}{d^{p(\theta-1)}}\,dx
\le
\gamma_p\Biggl\{&
\Bigl|\frac{p(\theta-1)+\beta(p-1)}{p}\Bigr|^p
\!\!\int_0^{+\infty}\! r^{-1}|g_\varepsilon(r)|^p\,dr\\
&+ c_1\!\int_0^{+\infty}\! |g_\varepsilon(r)|^{p-1}|g_\varepsilon'(r)|\,dr
+ c_2\!\int_0^{+\infty}\! r^{p-1}|g_\varepsilon'(r)|^p\,dr
\Biggr\}.
\end{split}
\]
By~\eqref{Integrali g_epsilon}, the last two integrals are uniformly bounded, while  
\[
\int_0^{+\infty}r^{-1}|g_\varepsilon(r)|^p\,dr
=-\ln(4\varepsilon^2)+\mathcal{O}(1).
\]
Hence,
\[
\int_{\mathbb{R}^N}
\Bigl|\nabla_{\L} u_\varepsilon \cdot 
\frac{\nabla_{\L} d}{|\nabla_{\L} d|}\Bigr|^{p}
\frac{1}{d^{p(\theta-1)}}\,dx
\le
\gamma_p\!\left(
-\Bigl|\frac{p(\theta-1)+\beta(p-1)}{p}\Bigr|^p\ln(4\varepsilon^2)
+\mathcal{O}(1)
\right),
\]
and similarly,
\[
\int_{\mathbb{R}^N}
\frac{|u_\varepsilon|^p}{d^{p\theta}}\,|\nabla_{\L} d|^p\,dx
=
\gamma_p\!\left(
-\ln(4\varepsilon^2)+\mathcal{O}(1)
\right).
\]
Dividing the two relations and letting \(\varepsilon\to0^+\) yields
\[
\lim_{\varepsilon\to0^+}
\frac{
\displaystyle
\int_{\mathbb{R}^N}
\Bigl|\nabla_{\L} u_\varepsilon \cdot 
\frac{\nabla_{\L} d}{|\nabla_{\L} d|}\Bigr|^{p}
\frac{1}{d^{p(\theta-1)}}\,dx
}{
\displaystyle
\int_{\mathbb{R}^N}
\frac{|u_\varepsilon|^p}{d^{p\theta}}\,|\nabla_{\L} d|^p\,dx
}
=
\Bigl|\frac{p(\theta-1)+\beta(p-1)}{p}\Bigr|^p.
\]
The same limit holds when \(|\nabla_{\L} u_\varepsilon|\) replaces  
\(|\nabla_{\L} u_\varepsilon \cdot \nabla_{\L} d / |\nabla_{\L} d||\),  
and thus the constant in both inequalities of~\eqref{INTRO: SV: Hardy Pesata} is sharp.
\end{proof}

\begin{proposition}[Sharpness of the weighted Rellich constant]\label{PROP: Rellich Sharpness}
Under the assumptions of Theorem~\ref{INTRO: SV: TH Rellich Pesata},
the constant in~\eqref{INTRO: SV: Rellich Pesata} is sharp.
\end{proposition}

\begin{proof}
Let 
\[
u_\varepsilon(x)
= d(x)^{\frac{\beta+p\theta+2p-2}{p}}\,g_\varepsilon(d(x)),
\qquad \varepsilon>0,
\]
where \(g_\varepsilon\) is the smooth cut-off function introduced in~\eqref{Funzioni Cut-Off}.  
A direct computation gives
\[
\begin{aligned}
\L u_\varepsilon
=&
-\frac{(2-\beta-p\theta-2p)(p\theta+2p-2-\beta(p-1))}{p^2}
\,d^{\frac{\beta+p\theta-2}{p}}\,g_\varepsilon(d)\,|\nabla_{\L} d|^2\\
&+
\Bigl(2\,\frac{\beta+p\theta+2p-2}{p}+1-\beta\Bigr)
\,d^{\frac{\beta+p\theta+p-2}{p}}\,g_\varepsilon'(d)\,|\nabla_{\L} d|^2
+
d^{\frac{\beta+p\theta+2p-2}{p}}\,g_\varepsilon''(d)\,|\nabla_{\L} d|^2.
\end{aligned}
\]

Using the coarea formula, the structural identity in 
Theorem~\ref{INTRO: SV: TH Rellich Pesata}, 
and Corollary~\ref{Stima di |a+b+c|^p}, we obtain
\[
\begin{aligned}
&\int_{\mathbb{R}^N} 
\frac{|\L u_\varepsilon|^p}{d^{p\theta}\,|\nabla_{\L} d|^{2(p-1)}}\,dx
\le
\lambda\Biggl\{
\left(
\frac{(2-\beta-p\theta-2p)\,(p\theta+2p-2-\beta(p-1))}{p^2}
\right)^p
\!\int_0^{+\infty} r^{-1}|g_\varepsilon(r)|^p\,dr\\
&\qquad
+\beta_1\!\int_0^{+\infty}\!|g_\varepsilon(r)|^{p-1}|g_\varepsilon'(r)|\,dr
+\beta_2\!\int_0^{+\infty}\!r\,|g_\varepsilon(r)|^{p-1}|g_\varepsilon''(r)|\,dr\\
&\qquad
+\beta_3\!\int_0^{+\infty}\!r^{p-1}|g_\varepsilon'(r)|^p\,dr
+\beta_4\!\int_0^{+\infty}\!r^p|g_\varepsilon'(r)|^{p-1}|g_\varepsilon''(r)|\,dr
+\beta_5\!\int_0^{+\infty}\!r^{2p-1}|g_\varepsilon''(r)|^p\,dr
\Biggr\},
\end{aligned}
\]
for suitable constants \(\beta_1,\dots,\beta_5>0\) depending only on \(p,\theta,\beta\).  

Hence, by~\eqref{Integrali g_epsilon},
\[
\int_{\mathbb{R}^N} 
\frac{|\L u_\varepsilon|^p}{d^{p\theta}\,|\nabla_{\L} d|^{2(p-1)}}\,dx
\le
\lambda\!\left\{
-\left(
\frac{(2-\beta-p\theta-2p)\,(p\theta+2p-2-\beta(p-1))}{p^2}
\right)^p\!\ln(4\varepsilon^2)
+\mathcal{O}(1)
\right\}.
\]
Similarly,
\[
\int_{\mathbb{R}^N}
\frac{|u_\varepsilon|^p}{d^{p(\theta+2)}}\,|\nabla_{\L} d|^2\,dx
=
\lambda\!\left(
-\ln(4\varepsilon^2)+\mathcal{O}(1)
\right).
\]
Dividing the two expressions and letting \(\varepsilon\to0^+\) gives
\[
\lim_{\varepsilon\to0^+}
\frac{
\displaystyle
\int_{\mathbb{R}^N} 
\frac{|\L u_\varepsilon|^p}{d^{p\theta}\,|\nabla_{\L} d|^{2(p-1)}}\,dx
}{
\displaystyle
\int_{\mathbb{R}^N}
\frac{|u_\varepsilon|^p}{d^{p(\theta+2)}}\,|\nabla_{\L} d|^2\,dx
}
=
\left(
\frac{(2-\beta-p\theta-2p)\,(p\theta+2p-2-\beta(p-1))}{p^2}
\right)^p,
\]
which establishes the sharpness of the constant in~\eqref{INTRO: SV: Rellich Pesata}.  
\end{proof}

\section{Applications and Model Cases}
\label{sezione: applicazioni}

In this section we apply the results established above to several classes of differential operators admitting a natural homogeneous structure. 
Before proceeding, we record two elementary identities concerning the \(p\)-Laplacian of powers and logarithms of a smooth function. 
These formulas will be used repeatedly throughout the examples.

\begin{lemma}
\label{SV: p-Laplaciano di d^beta}
Let \(D\subseteq \mathbb{R}^N\) be a domain, \(p\ge 2\), \(\beta\in\mathbb{R}\setminus\{0\}\), and \(f\in C^2(D)\) a positive function such that \(\L_p(f^\beta)=0\) in \(D\). 
Then
\[
\L_p(f)=\frac{(1-\beta)(p-1)}{f}\,|\nabla_{\L} f|^p 
\quad \text{in } D.
\]
\end{lemma}

\begin{proof}
For every smooth vector field \(\mathfrak{h}\), one has 
\(-\nabla_{\L}^*(f\mathfrak{h})=\nabla_{\L}f\cdot\mathfrak{h}-f\,\nabla_{\L}^*\mathfrak{h}\).
Taking 
\[
\mathfrak{h}
= |\beta|^{p-2}\beta\,
  f^{(\beta-1)(p-1)}
  |\nabla_{\L}f|^{p-2}\nabla_{\L}f,
\]
yields
\[
\L_p(f^\beta)
=|\beta|^{p-2}\beta
\bigl[
(\beta-1)(p-1)f^{(\beta-1)(p-1)-1}|\nabla_{\L}f|^p
+f^{(\beta-1)(p-1)}\L_p(f)
\bigr].
\]
Since \(f>0\) and \(\beta\neq0\), the prefactor is nonzero, and the condition
\(\L_p(f^\beta)=0\) implies
\[
\L_p(f)=\frac{(1-\beta)(p-1)}{f}\,|\nabla_{\L}f|^p.
\qedhere
\]
\end{proof}

\begin{lemma}
\label{SV: p-Laplaciano di ln}
Let \(D\subseteq \mathbb{R}^N\) be a domain, \(p\ge 2\), and \(f\in C^2(D)\) a positive function such that \(\L_p(-\ln f)=0\) in \(D\). 
Then
\[
\L_p(f)=\frac{(p-1)}{f}\,|\nabla_{\L}f|^p 
\quad \text{in } D.
\]
\end{lemma}

\begin{proof}
A direct computation shows that
\[
\L_p(-\ln f)
=\frac{(p-1)}{f^p}|\nabla_{\L}f|^p-\frac{1}{f^{p-1}}\L_p(f),
\]
and hence, if \(\L_p(-\ln f)=0\) in \(D\),
\[
\L_p(f)=\frac{(p-1)}{f}\,|\nabla_{\L}f|^p.
\qedhere
\]
\end{proof}

All the operators considered below admit a family of anisotropic dilations with respect to which they are homogeneous of degree two.  
Accordingly, throughout this section we assume that the function \(d\) introduced in the previous section satisfies this structure.  
That is, there exists a family of dilations \(\{\delta_\lambda\}_{\lambda>0}\),
\[
\delta_\lambda(x_1,\ldots,x_N)
=(\lambda^{\sigma_1}x_1,\ldots,\lambda^{\sigma_N}x_N),
\qquad
1\le \sigma_1\le\ldots\le\sigma_N,
\]
such that \(d\) and the horizontal gradient \(\nabla_{\L}=(X_1,\ldots,X_h)\) are homogeneous of degree one:
\[
d(\delta_\lambda(x))=\lambda d(x), 
\qquad
X_i(u\circ \delta_\lambda)=\lambda\,(X_i u)\circ \delta_\lambda, 
\quad i=1,\ldots,h.
\]
It follows that \(|\nabla_{\L}d|\) is homogeneous of degree zero.  
Consequently,
\[
\int_{\{d\le r\}}|\nabla_{\L}d|^p\,dx
=r^Q\int_{\{d\le1\}}|\nabla_{\L}d|^p\,dx
=\lambda_p\,r^Q,
\]
where \(Q=\sigma_1+\ldots+\sigma_N\) denotes the homogeneous dimension.  
Differentiating the coarea formula yields
\begin{equation}
\label{ES: Misura superficie Palla d}
\int_{\{d=r\}}\frac{|\nabla_{\L}d|^p}{|\nabla d|}\,d\mathcal{H}_{N-1}
=Q\,\lambda_p\,r^{Q-1}.
\end{equation}
This identity will be used to verify the integral conditions ensuring the sharpness of the constants in the subsequent examples.

\subsection{The Euclidean Laplacian}

When \(\L=\Delta\) and \(d(x)=|x|\), one has \(|\nabla d|=1\) and 
\[
\Delta_p d = \frac{N-1}{|x|} = \frac{(1-\beta)(p-1)}{d}\,|\nabla d|^p,
\]
so that the structural condition in Theorem~\ref{INTRO: SV: TH Hardy Pesata} is satisfied for every \(p\ge2\) with 
\[
\beta = \frac{p-N}{p-1}.
\]
All the assumptions of Theorems~\ref{INTRO: SV: TH Hardy Pesata} and~\ref{INTRO: SV: TH Rellich Pesata} are therefore trivially fulfilled with homogeneous dimension \(Q=N\).

As a consequence, for every \(u\in C_c^\infty(\mathbb{R}^N\setminus\{0\})\), \(p\ge2\), and \(\theta\in\mathbb{R}\),
\begin{equation}
\label{EUCL: Hardy}
\left|\frac{N-p\theta}{p}\right|^p
\int_{\mathbb{R}^N}\frac{|u|^p}{|x|^{p\theta}}\,dx
\le
\int_{\mathbb{R}^N}
\Bigl|\nabla u\cdot\frac{x}{|x|}\Bigr|^p
\frac{1}{|x|^{p(\theta-1)}}\,dx
\le
\int_{\mathbb{R}^N}\frac{|\nabla u|^p}{|x|^{p(\theta-1)}}\,dx.
\end{equation}
Moreover, whenever 
\((N-p(\theta+2))(p\theta+N(p-1))\ge0\),
\begin{equation}
\label{EUCL: Rellich}
\left(\frac{(N-p(\theta+2))(p\theta+N(p-1))}{p^2}\right)^p
\int_{\mathbb{R}^N}\frac{|u|^p}{|x|^{p(\theta+2)}}\,dx
\le
\int_{\mathbb{R}^N}\frac{|\Delta u|^p}{|x|^{p\theta}}\,dx.
\end{equation}
In both cases, the constants are sharp, and the corresponding (non-admissible) extremal profiles reduce, up to normalization, to the classical power functions
\[
u(x)=|x|^{\frac{p\theta-N}{p}}
\quad\text{and}\quad
u(x)=|x|^{\frac{p(\theta+2)-N}{p}},
\]
respectively.

\subsection{The Heisenberg--Greiner operator}

On \(\mathbb{R}^{2n+1}\) we consider the vector fields
\[
X_i = \frac{\partial}{\partial x_i} + 2\gamma\,y_i\,|z|^{2\gamma-2}\frac{\partial}{\partial t},
\qquad
Y_i = \frac{\partial}{\partial y_i} - 2\gamma\,x_i\,|z|^{2\gamma-2}\frac{\partial}{\partial t},
\quad i=1,\ldots,n,
\]
where \(z=(x,y)\in\mathbb{R}^{2n}\) and \(\gamma\ge1\).
The corresponding operator
\[
\L = \sum_{i=1}^n (X_i^2 + Y_i^2)
\]
is homogeneous of degree 2 with respect to the dilations
\[
\delta_\lambda(z,t) = (\lambda z, \lambda^{2\gamma} t), \qquad \lambda>0,
\]
which determine the homogeneous dimension \(Q = 2n + 2\gamma\). For \(\gamma = 1\), the operator reduces to the standard sub-Laplacian on the Heisenberg group,
while for integer values \(\gamma \ge 2\) the operator is a Greiner operator (see \cite{Gr}).

The associated homogeneous distance is
\[
\rho(z,t) = (|z|^{4\gamma} + t^2)^{1/(4\gamma)},
\qquad
|\nabla_{\L}\rho| = \frac{|z|^{2\gamma-1}}{\rho^{2\gamma-1}}.
\]
A \(p\)-harmonic function associated with \(\L_p\) is
\[
\Gamma_p(z,t) =
\begin{cases}
\rho^{\frac{p-Q}{p}}, & p\neq Q,\\[4pt]
-\ln\rho, & p=Q,
\end{cases}
\]
which satisfies \(\L_p\Gamma_p=0\) in \(\mathbb{R}^{2n+1}\setminus\{0\}\) (see \cite{ZN, DA0}).
By Lemmas~\ref{SV: p-Laplaciano di d^beta}--\ref{SV: p-Laplaciano di ln},
\[
\L_p\rho = \frac{(Q-1)(p-1)}{\rho}\,|\nabla_{\L}\rho|^p,
\]
so that the structural conditions of Theorems~\ref{INTRO: SV: TH Hardy Pesata}
and~\ref{INTRO: SV: TH Rellich Pesata} hold for all \(p\ge2\).

\begin{corollary}[Hardy inequality for the Heisenberg--Greiner operator]
\label{HG:Hardy}
Let \(p \ge 2\) and \(\theta \in \mathbb{R}\).
For every \(u \in C_c^\infty(\mathbb{R}^{2n+1}\setminus\{0\})\),
\[
\left|\frac{Q - p\theta}{p}\right|^p
\int_{\mathbb{R}^{2n+1}} 
\frac{|u|^p}{\rho^{p\theta}}\,|\nabla_{\L}\rho|^p\,dzdt
\le
\int_{\mathbb{R}^{2n+1}}
\left|
\nabla_{\L}u \cdot \frac{\nabla_{\L}\rho}{|\nabla_{\L}\rho|}
\right|^p
\frac{dzdt}{\rho^{p(\theta-1)}}
\le
\int_{\mathbb{R}^{2n+1}}
\frac{|\nabla_{\L}u|^p}{\rho^{p(\theta-1)}}\,dzdt.
\]
The constant on the left-hand side is optimal and never attained; equality holds,
up to scaling, for \(u = \rho^{\frac{p\theta - Q}{p}}\).
\end{corollary}

\begin{corollary}[Rellich inequality for the Heisenberg--Greiner operator]
\label{HG:Rellich}
Let \(p \ge 2\) and \(\theta \in \mathbb{R}\) satisfy
\[
(Q - p(\theta + 2))(p\theta + Q(p - 1)) \ge 0,
\qquad
Q > 2p(2\gamma - 1) + 4(1 - \gamma).
\]
Then, for all \(u \in C_c^\infty(\mathbb{R}^{2n+1}\setminus\{0\})\),
\[
\left(
\frac{(Q - p(\theta + 2))(p\theta + Q(p - 1))}{p^2}
\right)^p
\int_{\mathbb{R}^{2n+1}} 
\frac{|u|^p}{\rho^{p(\theta + 2)}}\,|\nabla_{\L}\rho|^2\,dzdt
\le
\int_{\mathbb{R}^{2n+1}}
\frac{|\L u|^p}{\rho^{p\theta}}\,|\nabla_{\L}\rho|^{-2(p - 1)}\,dzdt.
\]
The constant is optimal and equality occurs, up to a multiplicative factor, for
\(u = \rho^{\frac{p(\theta + 2) - Q}{p}}\).
\end{corollary}

The additional assumption on \(Q\) ensures the local integrability of the
right-hand side of the Rellich inequality.

\medskip
When \(\gamma = 1\), one recovers the sub-Laplacian on the Heisenberg group,
\[
\Delta_{\mathbb{H}^n}u
= \Delta_z u
+ 4|z|^2\frac{\partial^2 u}{\partial t^2}
+ 4\,\frac{\partial}{\partial t}(T u),
\qquad
T = \sum_{j=1}^n \left(y_j\frac{\partial}{\partial x_j}
- x_j\frac{\partial}{\partial y_j}\right),
\]
and the previous inequalities coincide with the classical Hardy and Rellich
inequalities on \(\mathbb{H}^n\) with sharp constants.

\subsection{The Baouendi--Grushin operator}

Let \(\gamma \ge 0\) and consider on \(\mathbb{R}^{n+k}\) the vector field
\[
\nabla_{\L} = \bigl(\nabla_x, (1+\gamma)|x|^{\gamma}\nabla_y\bigr),
\qquad
(x,y)\in\mathbb{R}^n\times\mathbb{R}^k,
\]
together with the corresponding operator
\[
\L = \Delta_x + (1+\gamma)^2 |x|^{2\gamma}\Delta_y.
\]
For \(k=0\) or \(\gamma=0\), this reduces to the classical Laplacian (and to the
\(p\)-Laplacian for the corresponding nonlinear operator \(\L_p\)).

The natural family of dilations associated with \(\L\) is
\[
\delta_\lambda(x,y) = (\lambda x, \lambda^{1+\gamma}y),
\qquad \lambda>0,
\]
with respect to which both \(\nabla_{\L}\) and \(\L\) are homogeneous of degree \(1\)
and \(2\), respectively.  
The corresponding homogeneous norm is
\[
\rho(x,y) = \bigl(|x|^{2(1+\gamma)} + |y|^2\bigr)^{1/(2+2\gamma)},
\qquad
|\nabla_{\L}\rho| = \frac{|x|^{\gamma}}{\rho^{\gamma}}.
\]
The homogeneous dimension is \(Q = n + (1+\gamma)k\).

Define
\[
\Gamma_p(x,y) =
\begin{cases}
\rho^{\frac{p-Q}{p}}, & p \neq Q,\\[4pt]
-\ln \rho, & p = Q,
\end{cases}
\]
which plays the role of the fundamental solution for \(\L_p\).  
A direct computation shows that \(\Gamma_p\) satisfies
\[
\L_p \Gamma_p = 0 \quad \text{in } \mathbb{R}^{n+k}\setminus\{0\}.
\]
By Lemmas~\ref{SV: p-Laplaciano di d^beta}--\ref{SV: p-Laplaciano di ln}, one has
\[
\L_p\rho = \frac{(Q - 1)(p - 1)}{\rho}\,|\nabla_{\L}\rho|^p,
\]
so that the structural assumptions of
Theorems~\ref{INTRO: SV: TH Hardy Pesata}--\ref{INTRO: SV: TH Rellich Pesata}
are satisfied for all \(p \ge 2\).

\begin{corollary}[Hardy inequality for the Baouendi--Grushin operator]
\label{BG:Hardy}
Let \(p \ge 2\) and \(\theta \in \mathbb{R}\).
Then, for every \(u \in C_c^\infty(\mathbb{R}^{n+k}\setminus\{0\})\),
\[
\left|\frac{Q - p\theta}{p}\right|^p
\int_{\mathbb{R}^{n+k}}
\frac{|u|^p}{\rho^{p\theta}}\,|\nabla_{\L}\rho|^p\,dxdy
\le
\int_{\mathbb{R}^{n+k}}
\left|
\nabla_{\L}u \cdot \frac{\nabla_{\L}\rho}{|\nabla_{\L}\rho|}
\right|^p
\frac{1}{\rho^{p(\theta-1)}}\,dxdy
\le
\int_{\mathbb{R}^{n+k}}
\frac{|\nabla_{\L}u|^p}{\rho^{p(\theta-1)}}\,dxdy.
\]
The constant on the left-hand side is optimal and never attained; equality holds,
up to scaling, for \(u = \rho^{\frac{p\theta - Q}{p}}\).
\end{corollary}

\begin{corollary}[Rellich inequality for the Baouendi--Grushin operator]
\label{BG:Rellich}
Let \(p \ge 2\) and \(\theta \in \mathbb{R}\) satisfy
\[
(Q - p(\theta + 2))(p\theta + Q(p - 1)) \ge 0\quad \text{and}\quad
Q>2\ga(p-1)+(1+\ga)k\,\,\text{ if }\,\, k\neq 0.
\]
Then, for all \(u \in C_c^\infty(\mathbb{R}^{n+k}\setminus\{0\})\),
\[
\left(
\frac{(Q - p(\theta + 2))(p\theta + Q(p - 1))}{p^2}
\right)^p
\int_{\mathbb{R}^{n+k}}
\frac{|u|^p}{\rho^{p(\theta + 2)}}\,|\nabla_{\L}\rho|^2\,dxdy
\le
\int_{\mathbb{R}^{n+k}}
\frac{|\L u|^p}{\rho^{p\theta}}\,|\nabla_{\L}\rho|^{-2(p - 1)}\,dxdy.
\]
The constant is optimal and equality occurs, up to a multiplicative factor, for
\(u = \rho^{\frac{p(\theta + 2) - Q}{p}}\).
\end{corollary}

When \(k = 0\) or \(\gamma = 0\), one recovers the Euclidean case.

\subsection{Homogeneous Carnot groups}
We briefly recall some basic facts on homogeneous Carnot groups and refer to \cite{BLU} for further background.

\begin{definition}
 We say that a Lie group on \(\rN\), 
\(G=(\rN,\circ)\), is a (homogeneous) Carnot group  if the following properties hold: 
\begin{enumerate}
\item \(\rN\) can be split as \(\rN=\erre^{N_1}\times\ldots\times\erre^{N_r}\),  and the dilation \(\delta_\la\colon \rN\to\rN\)
\[
\delta_\la(x)=\delta_\la(x^{(1)},\ldots,x^{(r)})
=
(\lambda x^{(1)},\lambda^2 x^{(2)},\ldots,\lambda^r x^{(r)})
\quad x^{(i)}\in\erre^{N_i},
\]
is an automorphism of the group \(G\) for every \(\la>0\).
\item  If \(N_1\) is as above, let \(X_1,\ldots,X_{N_1}\) 
be the left invariant vector fields on \(G\) such
that \(X_j(0)=\frac{\pa}{\pa x_j}\biggl|_0\) for \(j=1,\ldots, N_1.\) Then
\[
\mathrm{rank}\left(\mathrm{Lie}\{X_1,\ldots,X_{N_1}\}(x)\right)=N\quad\text{  for every } x\in\rN.
\]
\end{enumerate}
\end{definition}

We denote by
\[
\nabla_{\L} = (X_1,\ldots,X_{N_1}), 
\qquad 
\L = \sum_{i=1}^{N_1} X_i^2,
\]
the horizontal gradient and the canonical sub-Laplacian on \(G\).
For \(p \ge 2\), the \(p\)-sub-Laplacian is defined by
\[
\L_p u = \sum_{i=1}^{N_1} X_i\!\left(|\nabla_{\L}u|^{p-2} X_i u\right).
\]

Carnot groups are endowed with their Haar measure, which coincides with
Lebesgue measure on \(\mathbb{R}^N\).  
The homogeneous dimension associated with the dilations \(\delta_\lambda\) is
\[
Q = N_1 + 2N_2 + \cdots + r N_r,
\]
and it satisfies the scaling property
\[
|\delta_\lambda(E)| = \lambda^{Q} |E|
\qquad\text{for all measurable }E \subset \mathbb{R}^N.
\]
Moreover, the horizontal vector fields \(X_1,\ldots,X_{N_1}\) are
\(\delta_\lambda\)-homogeneous of degree \(1\).
From this point on we assume \(Q>3\); when \(Q \le 3\), Carnot groups reduce to the Euclidean setting and no genuinely
sub-Riemannian phenomena arise.

\medskip
Following \cite[Def.\ 5.1.1]{BLU}, we call a \emph{homogeneous norm} on \(G\) any continuous function 
\(\rho : G \to [0,+\infty)\) such that
\begin{itemize}
\item \(\rho(\delta_\lambda x) = \lambda\,\rho(x)\) for all \(\lambda>0\) and \(x \in G\);
\item \(\rho(x) > 0\) whenever \(x \ne 0\).
\end{itemize}

It is known (see, for instance, \cite[Ch.~5]{BLU}) that every Carnot group
admits a smooth homogeneous norm \(\rho_{2}\in C^{\infty}(G\setminus\{0\})\)
such that
\[
\Gamma(x)=\rho_{2}(x)^{\,2-Q}
\]
coincides, up to a multiplicative constant, with the fundamental solution of the sub-Laplacian \(-\L_{2}\).
Within this framework Balogh and Tyson introduced in \cite{BT} a notion of polarizability.
D’Ambrosio subsequently proved in \cite{DA0} that their notion is equivalent to the one recalled below.

\begin{definition}
A Carnot group \(G\) is called \emph{polarizable} if, for every \(p\ge 2\), the function
\[
\Gamma_{p}(x)=
\begin{cases}
\rho_{2}(x)^{\frac{p-Q}{p-1}}, & p\neq Q,\\[4pt]
-\ln \rho_{2}(x), & p=Q,
\end{cases}
\]
is \(\L_{p}\)-harmonic in \(G\setminus\{0\}\).
\end{definition}

Typical examples of polarizable Carnot groups are the Euclidean group
\((\mathbb{R}^{N},+)\), the Heisenberg group \(\mathbb{H}^{n}\), and H-type groups.

\medskip
By Lemmas~\ref{SV: p-Laplaciano di d^beta}–\ref{SV: p-Laplaciano di ln},
the homogeneous norm \(\rho_{2}\) associated with the sub-Laplacian \(\L_{2}\)
satisfies
\[
\L_{p}\rho_{2}
=
\frac{(Q-1)(p-1)}{\rho_{2}}\,
|\nabla_{\L}\rho_{2}|^{p}
\qquad\text{in }G\setminus\{0\}.
\]
This identity holds on every Carnot group when \(p=2\),
and for all \(p>2\) precisely when \(G\) is polarizable.

\begin{corollary}[Hardy inequality on Carnot groups]
\label{ES: TH: Hardy Carnot}
Let \(G\) be a Carnot group, and let \(p\ge2\) and \(\theta\in\mathbb{R}\).
Then every \(u\in C_{c}^{\infty}(G\setminus\{0\})\) satisfies:
\begin{itemize}
\item[\emph{(i)}] If \(p=2\),
\[
\left|\frac{Q-2\theta}{2}\right|^{2}
\int_{G}
\frac{|u|^{2}}{\rho_{2}^{2\theta}}\,|\nabla_{\L}\rho_{2}|^{2}\,dx
\le
\int_{G}
\left|
\nabla_{\L}u\cdot \frac{\nabla_{\L}\rho_{2}}{|\nabla_{\L}\rho_{2}|}
\right|^{2}
\frac{1}{\rho_{2}^{2(\theta-1)}}\,dx
\le
\int_{G}
\frac{|\nabla_{\L}u|^{2}}{\rho_{2}^{2(\theta-1)}}\,dx.
\]
The constant is sharp and equality occurs, up to a multiplicative factor, for 
\(u = \rho_{2}^{\frac{2\theta - Q}{2}}\).

\item[\emph{(ii)}] If \(p>2\) and \(G\) is polarizable,
\[
\left|\frac{Q-p\theta}{p}\right|^{p}
\int_{G}
\frac{|u|^{p}}{\rho_{2}^{p\theta}}\,|\nabla_{\L}\rho_{2}|^{p}\,dx
\le
\int_{G}
\left|
\nabla_{\L}u\cdot \frac{\nabla_{\L}\rho_{2}}{|\nabla_{\L}\rho_{2}|}
\right|^{p}
\frac{1}{\rho_{2}^{p(\theta-1)}}\,dx
\le
\int_{G}
\frac{|\nabla_{\L}u|^{p}}{\rho_{2}^{p(\theta-1)}}\,dx.
\]
The constant is sharp and equality occurs, up to a multiplicative factor, for 
\(u = \rho_{2}^{\frac{p\theta - Q}{p}}\).
\end{itemize}
\end{corollary}

\begin{corollary}[Rellich inequality on Carnot groups]
\label{ES: TH: Rellich Carnot}
Let \(G\) be a Carnot group, \(p\ge2\), \(\theta\in\mathbb{R}\) satisfy
\[
(Q-p(\theta+2))(p\theta+Q(p-1))\ge0,
\qquad
\rho_{2}^{-p\theta}\,|\nabla_{\L}\rho_{2}|^{-2(p-1)}
\in L^{1}_{\mathrm{loc}}(G).
\]
Then every \(u\in C_{c}^{\infty}(G\setminus\{0\})\) satisfies
\[
\left(
\frac{(Q-p(\theta+2))(p\theta+Q(p-1))}{p^{2}}
\right)^{p}
\int_{G}
\frac{|u|^{p}}{\rho_{2}^{p(\theta+2)}}\,|\nabla_{\L}\rho_{2}|^{2}\,dx
\le
\int_{G}
\frac{|\L u|^{p}}
{\rho_{2}^{p\theta}\,|\nabla_{\L}\rho_{2}|^{\,2(p-1)}}\,dx.
\]
The constant is sharp and equality occurs, up to a multiplicative factor, for 
\(u = \rho_{2}^{\frac{p(\theta + 2) - Q}{p}}\).
\end{corollary}

\nocite{*}
\bibliographystyle{plain}
\bibliography{Bibliografia_main}
\end{document}